\documentclass[a4paper,11pt,reqno]{amsart}
\usepackage{a4wide}
\usepackage{amsfonts,amsmath,amssymb}
\usepackage{url}

\newtheorem{conj}{Conjecture}
\theoremstyle{remark}

\begin{document}

\thispagestyle{plain}

\title{Asymptotics of generalised trinomial coefficients}

\author{Stephan Wagner}

\date{\today}

\begin{abstract}
It is shown how to obtain an asymptotic expansion of the generalised central trinomial coefficient $[x^n](x^2+b x+c)^n$ by means of singularity analysis, thus proving a conjecture of Zhi-Wei Sun.
\end{abstract}

\maketitle

In \cite{sun}, Zhi-Wei Sun proposes a number of conjectural formulas for multiples of $1/\pi$ involving the generalised central trinomial coefficient $T_n(b,c)$, which is defined by
$$T_n(b,c) = [x^n](x^2+bx+c)^n.$$
Besides the conjectural series, for example
$$\sum_{k=0}^{\infty} \frac{30k+7}{(-256)^k} \binom{2k}{k}^2 T_k(1,16) = \frac{24}{\pi},$$
he also presents conjectures regarding the asymptotic behaviour of $T_n(b,c)$ as $n \to \infty$, namely:

\begin{conj}[Sun {\cite[Conjecture 5.1]{sun}}]\label{sun}
For $b > 0$ and $c > 0$, we have
$$T_n(b,c) = \frac{(b+2\sqrt{c})^{n+1/2}}{2\sqrt[4]{c} \sqrt{n\pi}} \left( 1 + \frac{b-4\sqrt{c}}{16n \sqrt{c}} + O \left( \frac{1}{n^2} \right) \right)$$
as $n \to \infty$. If $c > 0$ and $b = 4\sqrt{c}$, then
$$\frac{T_n(b,c)}{\sqrt{c}^n} = \frac{3 \cdot 6^n}{\sqrt{6n\pi}} \left( 1 + \frac{1}{8n^2} + \frac{15}{64n^3} + \frac{21}{32n^4} + O \left( \frac{1}{n^5} \right) \right).$$
Finally, if $c < 0$ and $b \in \mathbb{R}$, then
$$\lim_{n \to \infty} \sqrt[n]{|T_n(b,c)|} = \sqrt{b^2-4c}.$$
\end{conj}

The aim of this little note is to show how the conjecture can be proven by a standard application of singularity analysis \cite{flaodl}.

\subsection*{The special cases $d = b^2-4c = 0$ and $b = 0$.} If the discriminant is $0$, then $T_n(b,c)$ essentially reduces to a central binomial coefficient:
$$T_n(b,c) = [x^n](x^2+bx+b^2/4)^n = [x^n](x+b/2)^{2n} = (b/2)^n \binom{2n}{n}.$$
In this case one can obtain an asymptotic expansion by means of Stirling's formula:
\begin{align*}
T_n(b,b^2/4) &= (b/2)^n \cdot \frac{(2n/e)^{2n}\cdot \sqrt{4\pi n}}{(n/e)^{2n}\cdot 2\pi n} \cdot \left( 1 + \frac{1}{24n} + \frac{1}{1152n^2} - \frac{139}{414720n^3} + O(n^{-4}) \right) \\
&\qquad \cdot \left( 1 + \frac{1}{12n} + \frac{1}{288n^2} - \frac{139}{51840n^3} + O(n^{-4}) \right)^{-2} \\
&= \frac{(2b)^n}{\sqrt{\pi n}} \left( 1 - \frac{1}{8n} + \frac{1}{128n^2} + \frac{5}{1024n^3} + O(n^{-4}) \right).
\end{align*}
Similarly, if $b = 0$, we obtain
$$T_n(0,c) = \begin{cases} c^{n/2} \binom{n}{n/2} & \text{$n$ even,} \\ 0 & \text{$n$ odd,} \end{cases}$$
and we can use Stirling's formula again.

\subsection*{The case $c > 0$.} Let us now assume that $b > 0$ (since $T_n(b,c) = (-1)^nT_n(-b,c)$, we can focus on this case), $c > 0$ and $d = b^2 - 4c \neq 0$. Then we can write
$$T_n(b,c) = d^{n/2} L_n(b/\sqrt{d}),$$
where $L_n$ is the $n$-th Legendre polynomial. It is thus sufficient to study the asymptotic behaviour of the Legendre polynomials. The Laplace-Heine formula states that
$$L_n(x) \sim \frac{(x+\sqrt{x^2-1})^{n+1/2}}{\sqrt{2n\pi}\sqrt[4]{x^2-1}}$$
as $n \to \infty$ if $-1 < x < 1$, which already yields the main term in Conjecture~\ref{sun} if $c > 0$. For our purposes, we mostly need the generating function
$$\sum_{n=0}^{\infty} L_n(x) t^n = \frac{1}{\sqrt{1-2xt+t^2}},$$
from which we get
\begin{align*}
F(t) &= \sum_{n=0}^{\infty} T_n(b,c) t^n = \sum_{n=0}^{\infty} L_n(b/\sqrt{d}) (\sqrt{d}t)^n \\
&= \frac{1}{\sqrt{1-2bt+dt^2}} = \frac{1}{\sqrt{1-2bt+(b^2-4c)t^2}}.
\end{align*}
The formula even remains valid when $d \leq 0$.
This function has two singularities at the zeros of the polynomial $1-2bt+(b^2-4c)t^2$. If $b > 0$ and $c > 0$, then these singularities are at $t_1 = 1/(b+2\sqrt{c})$ and at $t_2 = 1/(b-2\sqrt{c})$, and $t_1$ is closer to the origin. We now invoke singularity analysis (see \cite[Chapter VI]{flajolet} for a detailed explanation of this technique) to obtain the asymptotic behaviour of $T_n(b,c)$ from the expansion around the dominant singularity $t_1$:
$$F(t) = \frac{1}{2c^{1/4}} (t_1 - t)^{-1/2} - \frac{b^2-4c}{16c^{3/4}} (t_1-t)^{1/2} + \frac{3(b^2-4c)^2}{256c^{5/4}}(t_1-t)^{3/2} + \cdots$$
We can translate each term according to the rule
\begin{multline*}
C (1-t/t_1)^{-\alpha} \mapsto C \binom{n+\alpha-1}{n} t_1^{-n} \\
= \frac{Ct_1^{-n} n^{\alpha-1}}{\Gamma(\alpha)} \left( 1 + \frac{a(a-1)}{2n} + \frac{a(a-1)(a-2)(3a-1)}{24n^2} + O(n^{-3}) \right)
\end{multline*}
to obtain
$$T_n(b,c) = \frac{t_1^{-n-1/2}}{2c^{1/4}\sqrt{n\pi}} \left( 1 + \frac{b-4\sqrt{c}}{16\sqrt{c}n} + \frac{(3b-4\sqrt{c})^2}{512cn^2} + O(n^{-3}) \right).$$
This proves the first part of Conjecture~\ref{sun}, even with an additional term in the asymptotic expansion. By including further terms in the expansion around $t_1$ and in the expansion of the binomial coefficients $\binom{n+\alpha-1}{n}$, one can obtain even more precise asymptotic formulas. Let us illustrate this in the case that $b = 4\sqrt{c}$, when we get
$$F(t) = \frac{1}{\sqrt{1-2bt+3b^2t/4}} = \frac{1}{\sqrt{(1-bt/2)(1-3bt/2)}}.$$
The expansion around the dominant singularity $t_1 = 2/(3b)$ is, with $u = 1 - 3bt/2$,
$$F(t) = \sqrt{\frac32} \left( u^{-1/2} - \frac14 u^{1/2} + \frac{3}{32} u^{3/2} - \frac{5}{128} u^{5/2} + \frac{35}{2048} u^{7/2} - \frac{63}{8192} u^{9/2} + \cdots \right).$$
Moreover, a more precise asymptotic expansion of the Gamma function yields
\begin{multline*}
C \binom{n+\alpha-1}{n} t_1^{-n}
= \frac{Ct_1^{-n} n^{\alpha-1}}{\Gamma(\alpha)} \left( 1 + \frac{\alpha(\alpha-1)}{2n} + \frac{\alpha(\alpha-1)(\alpha-2)(3\alpha-1)}{24n^2} \right. \\
+ \frac{\alpha^2(\alpha-1)^2(\alpha-2)(\alpha-3)}{48n^3} +
\frac{\alpha(\alpha-1)(\alpha-2)(\alpha-3)(\alpha-4)(15\alpha^3-30\alpha^2+5\alpha+2)}{5760n^4} \\
\left. + \frac{\alpha^2(\alpha-1)^2(\alpha-2)(\alpha-3)(\alpha-4)(\alpha-5)(3\alpha^2-7\alpha-2)}{11520n^5} + O(n^{-6}) \right).
\end{multline*}
Hence the terms in the expansion of $F(t)$ translate as follows:
\begin{align*}
u^{-1/2} &\mapsto \frac{t_1^{-n}n^{-1/2}}{\Gamma(1/2)} \left( 1 - \frac{1}{8n} + \frac{1}{128n^2} + \frac{5}{1024n^3} - \frac{21}{32768n^4} - \frac{399}{262144n^5} + O(n^{-6}) \right), \\ 
u^{1/2} &\mapsto \frac{t_1^{-n}n^{-3/2}}{\Gamma(-1/2)} \left( 1 + \frac{3}{8n} + \frac{25}{128n^2} + \frac{105}{1024n^3} + \frac{1659}{32768n^4} + O(n^{-5}) \right), \\
u^{3/2} &\mapsto \frac{t_1^{-n}n^{-5/2}}{\Gamma(-3/2)} \left( 1 + \frac{15}{8n} + \frac{385}{128n^2} + \frac{4725}{1024n^3} + O(n^{-4}) \right), \\
u^{5/2} &\mapsto \frac{t_1^{-n}n^{-7/2}}{\Gamma(-5/2)} \left( 1 + \frac{35}{8n} + \frac{1785}{128n^2} + O(n^{-3}) \right), \\
u^{7/2} &\mapsto \frac{t_1^{-n}n^{-9/2}}{\Gamma(-7/2)} \left( 1 + \frac{63}{8n} + O(n^{-2}) \right), \\
u^{9/2} &\mapsto \frac{t_1^{-n}n^{-11/2}}{\Gamma(-9/2)} \left( 1 + O(n^{-1}) \right).
\end{align*}

Putting everything together, we arrive at
$$T_n(b,c) = \frac{3 \cdot 6^n}{\sqrt{6n\pi}} \left( 1 + \frac{1}{8n^2} + \frac{15}{64n^3} + \frac{21}{32n^4} + \frac{315}{128n^5} + O \left( \frac{1}{n^6} \right) \right),$$
which is the second part of Conjecture~\ref{sun}, even with one extra term.

\subsection*{The case $c < 0$.}
If $b > 0$ and $c < 0$, then the generating function $F(t)$ has two dominant singularities, since $t_1 = \frac{1}{b + 2i \sqrt{-c}}$ and $t_2 = \frac{1}{b - 2i \sqrt{-c}}$ have the same distance from the origin.
The singularity at $t_1$ is of the form
$$F(t) \sim \sqrt{\frac{b + 2i\sqrt{-c}}{4i\sqrt{-c}}} (1-t/t_1)^{-1/2},$$
and at $t_2$, it is
$$F(t) \sim \sqrt{\frac{b - 2i\sqrt{-c}}{-4i\sqrt{-c}}} (1-t/t_2)^{-1/2}.$$
Combining the contributions of the two, we obtain
\begin{align*}
T_n(b,c) &\sim \frac{1}{\Gamma(1/2)\sqrt{n}} \cdot \left( \sqrt{\frac{b + 2i\sqrt{-c}}{4i\sqrt{-c}}} \cdot t_1^{-n} + \sqrt{\frac{b - 2i\sqrt{-c}}{-4i\sqrt{-c}}} \cdot t_2^{-n} \right) \\
&= \frac{1}{2(-c)^{1/4}\sqrt{\pi n}} \cdot \left( \frac{1-i}{\sqrt{2}} (b+2i\sqrt{-c})^{n+1/2} + \frac{1+i}{\sqrt{2}} (b-2i\sqrt{-c})^{n+1/2} \right) \\
&= \frac{1}{(-c)^{1/4}\sqrt{\pi n}} (b^2-4c)^{n/2+1/4} \cos \left((n+1/2)\phi - \pi/4 \right),
\end{align*}
where $\phi$ is given by $e^{i\phi} = (b+2i\sqrt{-c})/\sqrt{b^2-4c}$. The third part of Conjecture~\ref{sun} follows immediately. Again, one could also get a further asymptotic expansion by taking more terms in the expansions around $t_1$ and $t_2$ into account.
 
\subsection*{Some final remarks} An alternative approach to these results would be the saddle point method (see Chapter VIII of \cite{flajolet}, cf. in particular Example VIII.11). Instead of the central trinomial coefficients, one can also treat other trinomial coefficients, see for instance \cite{katzenpanny1,katzenpanny2} for an application to random walks. Higher polynomial coefficients have of course been studied as well, see for example \cite[p.77]{comtet}.


\begin{thebibliography}{99}
 
\bibitem{comtet}
L.~Comtet,
\newblock {\em Advanced Combinatorics}.
\newblock D. Reidel Publishing Co., Dordrecht, 1974.

\bibitem{flaodl}
P.~Flajolet, A.~M. Odlyzko, 
\newblock Singularity analysis of generating functions,
\newblock {\em  SIAM J. Discrete Math.} 3 (1990), 216--240.

\bibitem{flajolet}
P.~Flajolet and R.~Sedgewick.
\newblock {\em Analytic combinatorics}.
\newblock Cambridge University Press, Cambridge, 2009.
\newblock Available online at \url{http://algo.inria.fr/flajolet/Publications/books.html}

\bibitem{katzenpanny1}
W.~Katzenbeisser, W.~Panny, 
\newblock Asymptotic results on the maximal deviation of simple random walks,
\newblock {\em Stochastic Process. Appl.} 18/2 (1984), 263--275.

\bibitem{katzenpanny2}
W.~Katzenbeisser, W.~Panny, 
\newblock Lattice path counting, simple random walk statistics, and randomization: an analytic approach. In \em{Advances in combinatorial methods and applications to probability and statistics}, 59–76, \em{Stat. Ind. Technol.}, Birkhäuser Boston, Boston, MA, 1997. 

\bibitem{sun}
Z.-W. Sun,
\newblock On sums related to central binomial and trinomial coefficients.
\newblock 2011.
\newblock Available online at \url{http://arxiv.org/abs/1101.0600v24}
\end{thebibliography}
\end{document}